\documentclass[11pt]{article}
\usepackage{amssymb,latexsym}
\usepackage{epsfig}
\usepackage{eufrak}
\usepackage{amsmath}
\usepackage{mathrsfs}
\usepackage{color}

\newtheorem{theorem}{Theorem}
\newtheorem{lemma}{Lemma}

\newtheorem{remark}{Remark}
\newtheorem{definition}{Definition}

\newcommand{\be}{\begin{equation}}
\newcommand{\ee}{\end{equation}}
\newcommand{\bee}{\begin{eqnarray*}}
\newcommand{\eee}{\end{eqnarray*}}
\newcommand{\bel}{\begin{eqnarray}}
\newcommand{\eel}{\end{eqnarray}}
\newcommand{\bec}{\begin{cases}}
\newcommand{\eec}{\end{cases}}
\newcommand{\bem}{\begin{bmatrix}}
\newcommand{\eem}{\end{bmatrix}}

\newcommand{\bed}{\begin{description}}
\newcommand{\eed}{\end{description}}
\newcommand{\bei}{\begin{itemize}}
\newcommand{\eei}{\end{itemize}}
\newcommand{\ben}{\begin{enumerate}}
\newcommand{\een}{\end{enumerate}}

\newcommand{\beL}{\begin{lemma}}
\newcommand{\eeL}{\end{lemma}}
\newcommand{\beT}{\begin{theorem}}
\newcommand{\eeT}{\end{theorem}}

\newcommand{\bpf}{\begin{pf}}
\newcommand{\epf}{\end{pf}}

\newcommand{\pfbox}{\hfill\mbox{$\Box$}}

\newenvironment{pf}{\paragraph*{Proof{\rm.}}}{\pfbox\bigskip}

\begin{document}

\title{{\bf Parallel Branch and Bound Algorithm for
Computing Maximal Structured Singular Value \thanks{This research
was supported in part by grants from AFOSR (F49620-94-1-0415),
ARO (DAAH04-96-1-0193), and LE{\cal Q}SF (DOD/LE{\cal
Q}SF(1996-99)-04).}}}
\author{Xinjia Chen and Kemin Zhou\\
Department of Electrical and Computer Engineering\\
Louisiana State University\\
Baton Rouge, LA 70803\\
chan@ece.lsu.edu \ \ kemin@ece.lsu.edu}

\date{January 1, 2002,  Revised on March 18, 2003\\
} \maketitle

\begin{abstract}

In this paper, we have developed a parallel branch and bound
algorithm which computes the maximal structured singular value
$\mu$ without tightly bounding $\mu$ for each frequency and thus
significantly reduce the computational complexity.

\end{abstract}

\begin{center}
\end{center}

\section{Introduction}

It is well known that the analysis of robust stability and
performance with structured uncertainty boils down to the problem
of computing the supremum of the structured singular value over
all frequency \cite{D,ZDG}. That is, $\mu_{max}:=\sup_{\omega \in
{\bf R}}\mu_{\bf \Delta}(M(j\omega))$ where $M(s)$ is the
transfer function of the generalized system and ${\bf \Delta}$ is
a set of block structured uncertainties.  Related to this problem
are the method proposed by Lawrence, Tits and Dooren \cite{Hel, LT} and the
approach established by  Ferreres and Biannic \cite{FB,FB2}.
These interesting techniques
can be applied to compute a $\mu$ upper bound over
a frequency interval without gridding of frequency.  For the precise
computation of the maximal structured singular value
$\mu_{max}$ (i.e., a tight lower bound is also expected in addition to an upper bound),
the conventional method is to grid a range
of frequency and compute the maximal $\mu$ among all the frequencies \cite{Bal}.
Since the exact computation is in general
impossible, $\mu$ is obtained for each frequency by tightly
bounding.  Sophisticated upper bounds and lower bounds have been
derived for example in \cite{B,D,FT,PD} and techniques such as
branch and bound \cite{NP}
 have been developed to refine the bounds.

It is noted that the existing techniques for computing the
maximal structured singular value $\mu_{max}$ lack of efficiency
because of the tedious frequency sweeping.  In this paper, we
investigate a smart frequency sweeping strategy. More
specifically, we apply branch and bound scheme to compute $\mu$
for $N > 1$ frequencies in parallel.  We introduce a powerful
``pruning '' mechanism.  That is, eliminate any branch with upper
bound smaller than $\frac{\hat{\mu}}{1-\epsilon}$ where
$\hat{\mu}$ is the maximum record of the lower bounds of all
branches ever generated and $\epsilon > 0$ is the tolerance. The
final $\hat{\mu}$ is returned as the maximal structured singular
value ${\mu}_{max}$. Since $\hat{\mu}$ is the maximum record of
the lower bounds obtained in all branches generated (no matter
belong to the same frequency or not), it will increase much
faster than its counterpart in the conventional frequency
sweeping algorithms.  Note that the raise of $\hat{\mu}$
results in a significant number
of branches to be pruned. Thus $\hat{\mu}$
 convergences quickly to the maximal
structured singular value ${\mu}_{max}$.

The paper is organized as follows. Section $2$ discusses existing
techniques for computing the maximal structured singular value.
Section $3$ presents our Parallel Branch and Bound Algorithm.  An
illustrative example is provided in Section $4$ and Section $5$
is the conclusion.

\section{Conventional Frequency Sweeping}

Consider an $M-\Delta$ set up as follows.

\begin{figure}[htb]
\centerline{\psfig{figure=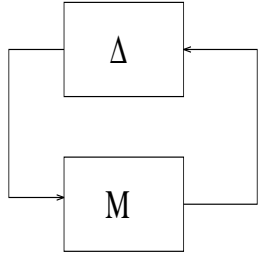,height=1.3in,width=1.3in}}
\caption{ Uncertain System.} \label{fig}
\end{figure}

Let ${\bf \Delta}$ be a set of block structured uncertainties. We
consider the computation of
\[
\mu_{max}:=\sup_{\omega \in {\bf R}}\mu_{\bf \Delta}(M(j\omega)).
\]
For notation simplicity, let $M(\omega)=M(j\omega)$.  Then
$\mu_{max}=\sup_{\omega \in {\bf R}}\mu_{\bf \Delta}(M(\omega))$.

In practice,  it is impossible to search $\mu_{max}$ over all
frequencies.  However, we can
 estimate $\mu_{max}$ as follows.

Choose a range of frequency $[a,\;b] \in {\bf R}$
and grid it as \be
\omega_j=a + \frac{(b-a)(j-1)}{NK-1},
\;\;\;j=1,\cdots, N K \label{grid} \ee where $N \geq 2$ and $K
\geq 1$ are integers (In practice, gridding is usually based on
the logarithmic scale.  However, in this paper, we use uniform
gridding for the simplicity of description.) Then an estimate for
$\mu_{max}$ can be defined as
\[
\tilde{\mu}_{max}:=\max_{j=1,\cdots, N K} \mu_{\bf
\Delta}(M(\omega_{j})).
\]

Define the (maximum positive real eigenvalue) function
$\bar{\lambda}_R: {\bf C}^{n \times n} \rightarrow {\bf R}$ as
\[
\bar{\lambda}_R (M):=\max \{\lambda:\;\lambda \;{\rm
is\;a\;positive\;real\;eigenvalue\;of}\;M\}
\]
with $\bar{\lambda}_R (M)=0$ if $M$ has no positive real
eigenvalues. Let ${\bf B \Delta}:= \{ \Delta \in {\bf \Delta}:
\bar{\sigma}(\Delta) \leq 1\}$. Then
\[
\mu_{\bf \Delta}(M)=\max_{\Delta \in {\bf B
\Delta}}\bar{\lambda}_R (M\Delta).
\]
Let $Q \subset {\bf B \Delta}$.  Define $\mu$ on a box \cite{NP}
\[
\mu(M,Q):= \max_{\Delta \in Q}\bar{\lambda}_R (M\Delta).
\]
There exists techniques in \cite{NP} for computing an upper bound
$UB(M,Q)$ and a lower bound $LB(M,Q)$ for $\mu(M,Q)$. Thus a
branch and bound scheme can be applied to compute $\mu_{\bf
\Delta}(M)$ with parameter space ${\bf B \Delta}$.

The conventional methods work essentially as follows.
For $j=1,\cdots,NK$, apply the following Algorithm $1$ to compute
an upper bound $UB^j$ and a lower bound $LB^j$ for $\mu_{\bf
\Delta}(M(\omega_j))$ such that $UB^j-LB^j \leq \varepsilon$.
Then $\tilde{\mu}_{max}$ satisfies
\[
\max_{j=1,\cdots,NK} LB^j \leq \tilde{\mu}_{max} \leq
\max_{j=1,\cdots,NK} UB^j
\]
where
\[
\max_{j=1,\cdots,NK} UB^j -\max_{j=1,\cdots,NK} LB^j \leq
\varepsilon.
\]

{\bf Algorithm $1$ --- Branch and Bound (\cite{NP})}

{\it Initialize}  Let ${\cal U}_{j}=\{Q_k\}={\bf B\Delta}$.

{\it Let}
\[
UB^j= \max_{k} UB(M(\omega_j),Q_k),
\]
\be LB^j= \max_{k} LB(M(\omega_j),Q_k). \label{global} \ee

{\it while} $UB^j-LB^j > \varepsilon$

\begin{itemize}

\item Choose $Q$ to be any element of ${\cal U}_{j}$
with $UB(M(\omega_{j}),Q)=UB^j$.

\item Partition $Q$ into $Q_a$ and $Q_b$ by
 bisecting along one of its longest edges.

\item Add $Q_a$ and $Q_b$ into ${\cal U}_{j}$.
Remove $Q$ from ${\cal U}_{j}$.

\item Remove from ${\cal U}_{j}$ any $Q$ with
\be UB(M(\omega_{j}),Q) < LB^j. \label{con2} \ee

\end{itemize}

{\it endwhile}

The most important mechanism of Algorithm $1$ is ``pruning''
\cite{NP}.  That is, any element of ${\cal U}_{j}$ for which
~(\ref{con2}) is satisfied will never again be partitioned and
need not be considered further. We call inequality ~(\ref{con2})
as the ``pruning condition''.

We can see that existing techniques for computing
$\tilde{\mu}_{max}$ employ branch and bound techniques for each
frequency independently.  In particular, the pruning process for
one frequency is independent of another. $\mu_{\bf \Delta}(M)$ is
bounded tightly for each frequency. Note that we usually need to
evaluate $\mu_{\bf \Delta}(M)$ for many
  frequencies in order to obtain a reasonably good estimate of
 the maximal structured singular value $\mu_{max}$.  Thus
 the overall computation is still a heavy burden,
 even the computation of $\mu_{\bf \Delta}(M)$ for each frequency is
 very efficient.

Thus for the sake of efficiency, there is a strong motivation to
conceive
 a smart frequency sweeping strategy.  More specifically,
 we would raise the following question,

{\it Is it possible to obtain the maximal structured singular
value $\mu_{max}$ without tightly bounding $\mu_{\bf
\Delta}(M(\omega_{j}))$ for each frequency $\omega_{j}$?}

The following section is devoted to answering this question.

\section{Parallel Branch and Bound Algorithm}
It is fair to compare the performance of different algorithms on
the same set of frequencies.  Therefore, we consider again
frequencies $\omega_j,\;\;j=1,\cdots,NK$ defined by ~(\ref{grid})
and relabel them as
\[
\omega_{ij}:= a + \frac{(b-a)[K(i-1) + (j-1)]}{NK-1}
,\;\;\;\;i=1,\cdots,N,\;\;j=1,\cdots,K.
\]
Now we are in a good position to present our Parallel Branch and
Bound Algorithm as follows.

{\bf Algorithm $2$ --- Parallel Branch and Bound Algorithm}
\begin{itemize}

\item Step $1$: Initialize. Set $j=1$. Set $\hat{\mu}=0$.
Set tolerance $\epsilon > 0$.  Set maximal iteration number $IT$.

\item Step $2$: Update $\hat{\mu}$ and record the number of
iterations $r(j)$ for frequency $\omega_{ij}$ by the following
steps.

\begin{itemize}

\item Step $2$--$1$:
Let ${\cal U}_{ij}=\{Q_k\}={\bf B\Delta},\;\; i=1,\cdots,N$. Set
$r=1$.

\item Step $2$--$2$: If $r=IT+1$ or ${\cal U}_{ij}$ is empty for any
$i \in \{1,\cdots,N\}$ then record $r(j)=r$ and go to Step $3$,
else do the following for all $i$ such that ${\cal U}_{ij}$ is
not empty.

\begin{itemize}

\item Choose $Q$ to be any element of ${\cal U}_{ij}$
with
\[
UB(M(\omega_{ij}),Q)= \max_{Q_k \in {\cal U}_{ij}}
\;UB(M(\omega_{ij}),Q_k).
\]

\item Partition $Q$ into $Q_a$ and $Q_b$ by
 bisecting along one of its longest edges.

\item Add $Q_a$ and $Q_b$ into ${\cal U}_{ij}$.
Remove $Q$ from ${\cal U}_{ij}$.

\item Update
\be \hat{\mu}=\max \{\hat{\mu},\;
LB(M(\omega_{ij}),Q_a),\;LB(M(\omega_{ij}),Q_b)\}. \label{update}
\ee

\item Remove from ${\cal U}_{ij}$ any $Q$ with
\be UB(M(\omega_{ij}),Q) < \frac{\hat{\mu}} {1-\epsilon}.
\label{con} \ee

\end{itemize}

\item Step $2$--$3$: Set $r=r+1$ and go to Step $2$--$2$.

\end{itemize}

\item Step $3$: If $j=K$ then STOP, else set $j=j+1$ and go to Step $2$.

\end{itemize}

In Algorithm $2$, $N$ branches of frequency sweeping are performed
in parallel with starting frequencies
$\omega_{i1},\;\;i=1,\cdots,N$ and step size
$\frac{b-a}{NK-1}$.  Also, a branch and bound
scheme is applied to compute $\mu$ for $N$ frequencies in
parallel.  Any branch with upper bound smaller than
$\frac{\hat{\mu}}{1-\epsilon}$ will be pruned, where $\hat{\mu}$
is the maximum record of the lower bounds of all branches ever
generated. The final $\hat{\mu}$ is returned as the maximal
structured singular value $\mu_{max}$.  Algorithm $2$ is visulized
in the following Figure ~\ref{fig_18}.

\begin{figure}[htb]
\centerline{\psfig{figure=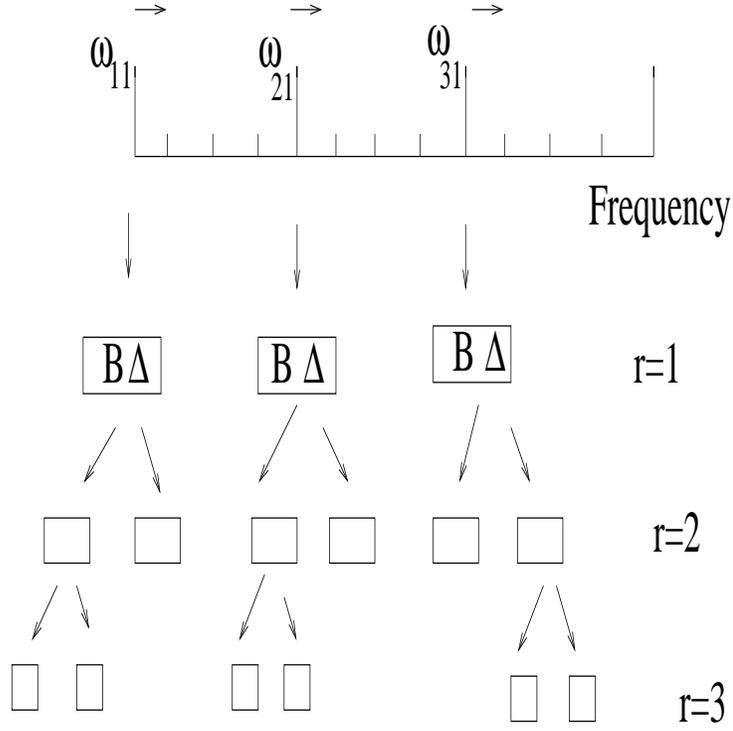,height=3.8in  ,width=3.8in  }}
\caption{A Picture of Parallel Branch and Bound Algorithm.
$N=3,\;\;K=4$.
 } \label{fig_18}
\end{figure}

\begin{remark}
Note that Algorithm $2$ provides a substantial improvement on
efficiency than conventional methods in computing the maximal
structured singular value. This can be explained by the
significant relaxation in the ``pruning condition'' of Algorithm
$2$. To see the difference of the two ``pruning conditions'', we
can compare the right hand sides of inequalities ~(\ref{con}) and
~(\ref{con2}). By ~(\ref{update}) and ~(\ref{global}), we can see
that $\frac{\hat{\mu}} {1-\epsilon}$ can be much larger than
$LB^j$.  This is because $\hat{\mu}$ is the maximum record of the
lower bounds obtained in all branches of all frequencies
evaluated and being evaluated, while $LB^j$ is only the maximum
record of the lower bounds obtained in branches of the frequency
being evaluated. Moreover, $\hat{\mu}$ is enlarged to
$\frac{\hat{\mu}} {1-\epsilon}$ in the ``pruning condition''
~(\ref{con}) and hence the ``pruning'' process is further
facilitated. The significant relaxation of the `pruning
condition'' leads to a substantial decrease of the number of
total subdomains needed to be evaluated. Therefore, our algorithm
is much more efficient than those previously available to control
engineers.

\end{remark}

\begin{remark}
It is important to note that Algorithm $2$ involves only one CPU
processor. It is fundamentally different from the parallel
algorithms which involves more than one CPU processors.
\end{remark}

\begin{remark}
A substantial amount of computation can be saved by the following
mechanisms.  First, further
 computation of the lower
bound on a domain is not needed once it is determined that the
lower bound is smaller than the existing global lower bound.
This can be seen from equation ~(\ref{update}).  Second, the
computation of the upper bound should be terminated once
condition ~(\ref{con}) is satisfied.   The idea of these two
mechanisms is to avoid as much as possible the tightly
computation of the lower bound and the upper bound.
\end{remark}

In addition to the novel frequency sweeping strategy, another
character of Algorithm $2$ is that there is no tolerance criteria
directly forced on the final result, however, the final result
falls into tolerance automatically.

\begin{theorem} \label{t4}
Suppose that the maximal iteration number $IT < \infty$ and that
Algorithm $2$ stops with $r(j) \leq IT,\;\;j=1,\cdots,K$. Then
the final $\hat{\mu}$ satisfies
\[
0 \leq \frac{\tilde{\mu}_{max}-\hat{\mu}} {\tilde{\mu}_{max}}
< \epsilon.
\]
\end{theorem}

\begin{pf} Since $\hat{\mu}$ is the maximal record of the lower bounds,
we have $\tilde{\mu}_{max}-\hat{\mu} \geq 0$.  We only need to
show that $\frac{\tilde{\mu}_{max}-\hat{\mu}}
{\tilde{\mu}_{max}}    < \epsilon$. By the assumption that
Algorithm $2$ stops with $r(j) \leq IT,\;\;j=1,\cdots,K$, we know
that all subdomains ever generated are finally removed because
the ``pruning condition'' ~(\ref{con}) is satisfied.  Note that
there exists a subdomains $Q_{ij}$ for frequency $\omega_{ij}$
such that $\mu(M(\omega_{ij}), Q_{ij})=\tilde{\mu}_{max}$. Let
$\hat{\mu}=\bar{\mu}$ when $Q_{ij}$ is removed. Then $
\tilde{\mu}_{max} \leq UB(M(\omega_{ij}),Q_{ij}) <
\frac{\bar{\mu}} {1-\epsilon}$. Note that $\hat{\mu}$ is
nondecreasing thus the final $\hat{\mu} \geq \bar{\mu}$.  It
follows that
\[
\tilde{\mu}_{max} < \frac{\hat{\mu}} {1-\epsilon} \;\;\;
\Longrightarrow \;\;\; \frac{\tilde{\mu}_{max}-\hat{\mu}}
{\tilde{\mu}_{max}}    < \epsilon.
\]
The proof is thus completed.

\end{pf}

Note that one important concern of an algorithm is convergence.
It is usually desirable that, given any tolerance $\epsilon > 0$,
an algorithm stops and returns the result within tolerance in a
finite number of iterations.
 Obviously, the convergence requirement imposes
 condition of the quality of bounds.

\begin{definition} \label{d1}
The upper bound $UB(M,.)$ and lower bound $LB(M,.)$ are said to
be continuous if
\[
\lim_{d(Q) \rightarrow 0} UB(M,Q) - LB(M,Q) =0
\]
where $d(Q): =\max_{q,\;q^{'} \in Q} ||q-q^{'}||$ with $Q
\subseteq {\bf B\Delta}$.
\end{definition}

\begin{theorem} \label{t5}
Suppose that all the upper bounds and lower bounds are continuous
and that at least one nonzero lower bound appears after a finite
number of iterations.  Let the maximal iteration number $IT = \infty$.  Then, for
arbitrary tolerance $\epsilon > 0$, Algorithm $2$ stops with a
finite number of domain partitions for each $j$, i.e., $r(j) <
\infty,\;\;j=1,\cdots,K$. Moreover, the final $\hat{\mu}$
satisfies
\[
0 \leq \frac{\tilde{\mu}_{max}-\hat{\mu}} {\tilde{\mu}_{max}}
< \epsilon.
\]
\end{theorem}

\begin{pf}
Suppose that Algorithm $2$ does not stop with a finite number of
domain partitions for each $j$. Then $\exists \;\omega_{ij}$ and
an infinite sequence of nested subdomains $\{Q_{r}^{ij}\}$
associated with frequency $\omega_{ij}$ such that $ Q_{1}^{ij}
\supset Q_{2}^{ij}  \supset \cdots  \supset Q_{r}^{ij} \supset
\cdots$. Note that by the assumption $\exists r_0 < \infty, \mu_0
> 0$ such that $\hat{\mu} \geq  \mu_{0},\;\;\forall r > r_0$. By
the continuity, $\exists r_1$ such that
\[
UB(M(\omega_{ij}),Q_{r}^{ij}) - LB(M(\omega_{ij}),Q_{r}^{ij}) <
\frac{\epsilon}{1-\epsilon} \mu_{0}, \;\;\forall r > r_1.
\]
Let $r_2=\max\{r_0,r_1\}+1$.  Then
\[
UB(M(\omega_{ij}),Q_{r_2}^{ij}) - LB(M(\omega_{ij}),Q_{r_2}^{ij})
< \frac{\epsilon}{1-\epsilon} \mu_{0}.
\]
Thus
\[
UB(M(\omega_{ij}),Q_{r_2+1}^{ij})-\hat{\mu} <
\frac{\epsilon}{1-\epsilon}\hat{\mu} \;\;\;\Longrightarrow
\;\;\;UB(M(\omega_{ij}),Q_{r_2+1}^{ij}) <
\frac{\hat{\mu}}{1-\epsilon}
\]
which implies that $Q_{r_2+1}^{ij}$ is removed. This is a
contradiction. Therefore Algorithm $2$ stops with a finite number
of domain partitions for each $j$ and hence by the same argument
of Theorem ~\ref{t4}
\[
0 \leq \frac{\tilde{\mu}_{max}-\hat{\mu}} {\tilde{\mu}_{max}}
< \epsilon.
\]
The proof is thus completed.
\end{pf}

\section{An Illustrative Example}

Consider an $M-\Delta$ set up as shown in Figure ~\ref{fig} where
$\Delta = {\rm diag} (\delta_1,\delta_2) \in {\bf R}^{2 \times 2}$ and
$M(s)=C(sI-A)^{-1}B$ with
\[
A=\left[\begin{array}{cccc}
-1 & -10 & -1 & 10\\
-0.5 & -1 & 1 & 0.5\\
0.5 & -4 & -1 & -10\\
-10 & 0.5 & 0 &
-2.5\end{array}\right],\;\;\;B=\left[\begin{array}{cc}
1 & 0\\
0 & 0\\
0 & 0\\
0 & 1\end{array}\right],\;\;\;C=\left[\begin{array}{cccc}
-0.5 & 0 & 0 & 0\\
0 & 0 & 0 & -1.5\\
\end{array}\right].
\]
To compute the supremum of $\mu$, we uniformly grid frequency
interval $[0.01 ,\;15.01]$ and obtain $1,500$ grid frequencies as
\[
\omega_j=0.01 + \frac{(15.01 - 0.01) (j-1)}{1500 -1}, \;\;j=1,\cdots,1500.
\]
In Algorithm $2$, we choose the relative error $\epsilon=0.01$
and $N=30,\;\;K=50$. The $1,500$ frequencies are regrouped as
\[
\omega_{ij}= 0.01 + \frac{(15.01 - 0.01) [ 50(i-1) + (j-1)] }{1500 -1},\;\;i=1,\cdots,30;\;\;j=1,\cdots,50.
\]
The execution of Algorithm 2 is terminated at $r = 1$ with $\hat{\mu}=0.8424$ achieved at frequency
$\omega_{ij}=9.1661$ where $i = 19, \;j = 16$.  It is observed that, for any frequency, no
partition is performed for the original domain ${\bf B \Delta} =[-1,1]\times[-1,1]$.
This is because $\frac{\hat{\mu}}{1 - \epsilon}$
is greater than the upper bounds of singular values for ${\bf B \Delta}$ at other frequencies.
The bounds of singular values for ${\bf B \Delta}$ are shown by Figures 3-4.
It can be seen from these figures that the upper bounds and lower bounds of  singular values
are far apart for most of the frequencies.  To compute the maximal singular value
using the conventional branch and bound method,
substantial computational effort will be wasted on reducing the gap between
the upper bounds and lower bounds of  singular values for most of the frequencies.
This example demonstrates that branch and bound should not be applied extensively for any fixed frequency.
Quiet contrary, it should be employed in parallel and in a cooperative manner.
This spirit has been reflected in Algorithm $2$.

\begin{figure}[htb]
\centerline{\psfig{figure=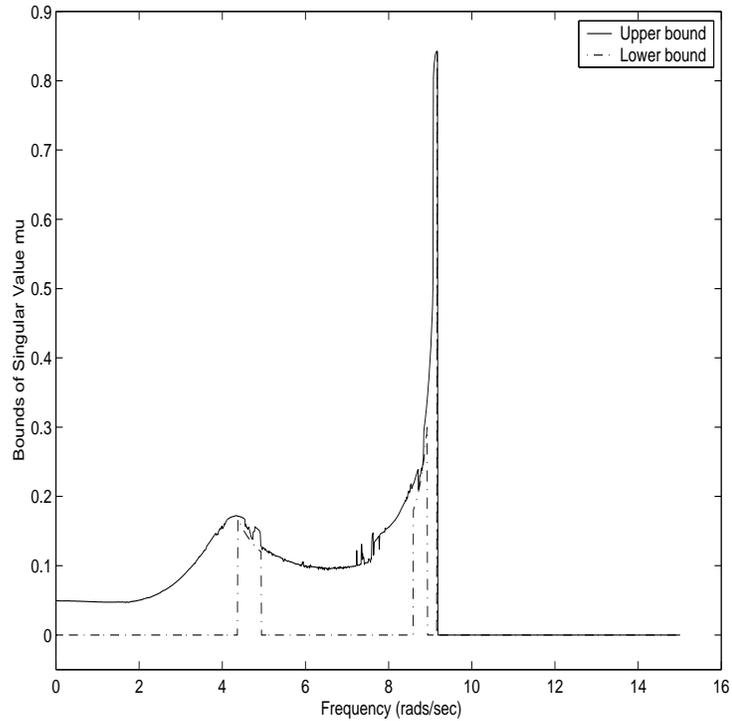,height=3.8in  ,width=3.8in
}} \caption{Bounds of Singular Values}
 \label{fig_10}
\end{figure}

\begin{figure}[htb]
\centerline{\psfig{figure=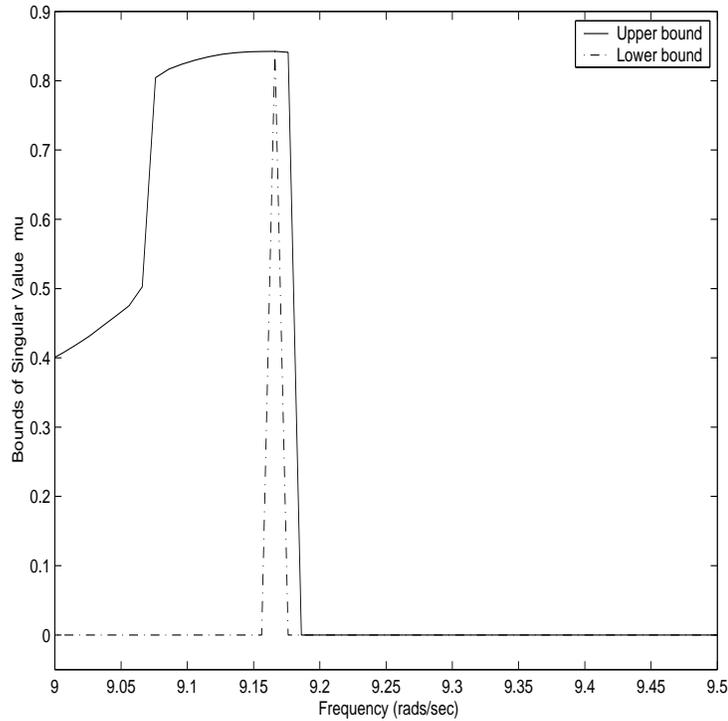,height=3.8in  ,width=3.8in
}} \caption{Bounds of Singular Values}
 \label{fig_12}
\end{figure}

\section{Conclusion}

Efficient computation of the maximal structured singular value is
of fundamental importance in robustness analysis and robust
synthesis with structured uncertainty. Motivated by this, we have
developed a parallel branch and bound algorithm for computing the
maximal structured singular value, which significantly reduce the
computational complexity.

\end{document}